\documentclass[12pt]{article}
\usepackage{mathtools} 
\usepackage{amsmath, amssymb} 
\usepackage{amsthm}
\usepackage{enumerate}  % lets you do \begin{enumerate}[(a)] and such
\usepackage{url} % formats url's nicely using \url{}
\usepackage{authblk}

\usepackage{tikz}
\usepackage{tkz-graph}
\usetikzlibrary{shapes}
\usetikzlibrary{arrows}
\usetikzlibrary{decorations.markings}

\usepackage{graphicx}
\usepackage{caption,subcaption}

\newcommand\cx{{\mathbb C}}% complexes

\newcommand\ints{{\mathbb Z}}
\newcommand\re{{\mathbb R}}%reals
\newcommand\rats{{\mathbb Q}}

\DeclarePairedDelimiter\abs{\lvert}{\rvert}%
\DeclarePairedDelimiter\norm{\lVert}{\rVert}%

\makeatletter
\let\oldabs\abs
\def\abs{\@ifstar{\oldabs}{\oldabs*}}
\let\oldnorm\norm
\def\norm{\@ifstar{\oldnorm}{\oldnorm*}}
\makeatother

\newcommand\sbs{\subseteq}

\newcommand\seq[4]{#1_{#2},#1_{#3},\ldots,#1_{#4}}

\newtheoremstyle{plainsl}%
	{\topsep}
	{\topsep}
	{\slshape} % only non-default setting
	{}
	{\normalfont\bfseries}
	{.}
	{ }
	{}

\swapnumbers

{\theoremstyle{plainsl}
\newtheorem{theorem}{Theorem}[section]
\newtheorem{lemma}[theorem]{Lemma}
\newtheorem{corollary}[theorem]{Corollary}

}
{\theoremstyle{remark}

}

\renewcommand\proof{\noindent\textsl{Proof. }}
\newcommand\sqr[2]{{\vbox{\hrule height.#2pt
    \hbox{\vrule width.#2pt height#1pt \kern#1pt
        \vrule width.#2pt}\hrule height.#2pt}}}
\renewcommand\qed{%
	\ifmmode\eqno\sqr53
	\else\nolinebreak\ \hfill\sqr53\medbreak\fi}

\DeclareMathOperator{\wt}{wt}

\newcommand\ip[2]{\left\langle#1,#2\right\rangle}
\newcommand\one{{\bf1}}
\newcommand\zero{{\bf0}}
\usepackage{blkarray}

\title{Simple Quantum Coins Enable Pretty Good State Transfer on Every Hypercube}
\author{Hanmeng Zhan}
\date{}

\affil{Computer Science Department\\ Worcester Polytechnic Institute, Worcester, MA, USA\\\texttt{hzhan@wpi.edu}}
\begin{document}
\maketitle

\begin{abstract}
We consider pretty good state transfer in coined quantum walks between antipodal vertices on the hypercube $Q_d$. When $d$ is a prime, this was proven to occur in the arc-reversal walk with Grover coins. We extend this result by constructing weighted Grover coins that enable pretty good state transfer on every $Q_d$. Our coins are real, and require modification of the weight on only one arc per vertex. We also generalize our approach and establish a sufficient condition for  pretty good state transfer to occur on other graphs.
\end{abstract}

\section{Introduction}
A number of quantum algorithms can be seen as discrete-time quantum walks on graphs with nice transport properties. For example, in a search algorithm (e.g. \cite{Grover1996} and \cite{Ambainis2001}), the quantum walk starts with a uniform superposition of all arcs of a graph and gets very close to a state that concetrates on a vertex. Properties like this motivated the study of various types of quantum state transfer, including perfect state transfer \cite{Lovett2009,Kendon2011,Kurzynski2011,Barr2014,Yalcnkaya2014,Stefanak2016,Stefanak2017,Zhan2019,Kubota2021,Guo2024a,Chen2024a,Kubota2024,Bhakta2024a}, periodicity \cite{Yoshie2019,Ito2020,Yoshie2023,Chen2024,Bhakta2024}, pretty good state transfer \cite{Chan2023}, and more recently peak state transfer \cite{Guo2024}. 

We are particularly interested in state transfer between antipodal vertices on the hypercube $Q_d$, where the quantum walk starts in a superposition of the outgoing arcs of a vertex and approaches a superposition of the outgoing arcs of the antipodal vertex. It was shown by \cite{Chan2023} that for $d$ prime, this transfer can be \textsl{pretty good}---for any $\epsilon>0$, there is a time at which the system is $\epsilon$-close to the target state---and we may achieve this in an arc-reversal walk with Grover coins. However, the same strategy does not work for all hypercubes: on $Q_4$,  the arc-reversal walk with Grover coins is periodic, and the maximum transfer probability between the antipodal vertices is only $9/16$. 

One way to boost the transfer probability is to mark the initial and target vertices \cite{Stefanak2023}, an idea inspired by quantum search algorithms. Alternatively, we may change the coins to obtain a quantum walk that satisfies the algebraic conditions for pretty good state transfer. In this paper, we show a slight modification of the Grover coin achieves our goal: for every $d$, there is a weighted Grover coin relative to which the arc-reversal walk on $Q_d$ exhibits pretty good state transfer between the antipodal vertices. Our work  relies heavily on the theory developed for pretty good state transfer in a class of quantum walks \cite{Chan2023}, where the coin operator is a reflection about certain complex-weighted tail-arc incidence matrix of the graph. 

\section{The quantum walk with respect to $(X,W)$}
We briefly describe the quantum walk introduced by \cite{Chan2023}. Given a connected graph $X$, and a complex-weighted adjacency matrix $W$ of $X$, where $W_{ab}\ne 0$ if and only if $a$ is adjacent to $b$, \textsl{the quantum walk with respect to $(X,W)$} takes place on the \textsl{arcs} of $X$:
\[\{(u,v): u, v\in V(X), \{u,v\}\in E(X)\},\]
where the shift operator reverses each arc, and the coin operator is determined by the entries of $W$. More specifically, let $R$ be the permutation matrix that swaps arc $(u,v)$ with arc $(v,u)$, and let $N_t$ be the \textsl{weighted arc-tail incidence matrix}, where 
\[(N_t)_{u, (a, b)} = \begin{cases}
	\dfrac{w_{ab}}{\sqrt{(WW^*)_{aa}}}, & \text{if }u=a,\\
	0,\quad & \text{if } u\ne a.
\end{cases}.\]
Then the transition matrix is the product
\[U=R(2N_t^*N_t-I).\]
In particular, if $W$ assigns the same weight to each arc, then $2N_t^*N_t-I$ is coin operator that assigns the Grover coin $\frac{2}{\deg(u)}J-I$ to the vertex $u$.

Let 
\[H=(I\circ WW^*)^{1/2} (W\circ W^*) (I\circ WW^*)^{1/2}.\]
It was shown in \cite{Chan2023} that
\[H = N_tRN_t^*,\]
and the spectral decomposition of $U$ is related to the spectral decomposition of $H$. For this reason, we will refer to $H$ as the \textsl{Hermitian adjacency matrix of $X$ associated with $W$}.

\begin{theorem}\label{thm:eprojs}\cite[Theorem 3.3]{Chan2023}
	Let $U$ be the transition matrix with respect to $(X,W)$, and let $H$ be the Hermitian adjacency matrix of $X$ associated with $W$. The eigenvalues of $U$ are $1$, $-1$ and $e^{\pm i\theta}$, where $\theta=\arccos\lambda$ for some eigenvalue $\lambda\in(-1,1)$ of $H$.  Moreover, if $F_{\theta}$ denotes the $e^{i\theta}$-eigenprojection of $U$, and $E_{\lambda}$ denotes the $\lambda$-eigenprojection of $H$, then $F_{\theta}$ and $E_{\lambda}$ are related in the following way.
	\begin{enumerate}[(i)]
		\item If $\theta=0$, then $\lambda=1$ and
		\[ N_t F_0 N_t^* = E_1.\]
		\item If $\theta=\pi$, then $\lambda=-1$ and
		\[N_t F_{\pi}N_t^* = E_{-1}.\]
		\item If $\theta \in(-\pi,0)\cup(0,\pi)$, then
		\[N_t F_{\theta}N_t^* = \frac{1}{2}E_{\lambda}.\]
	\end{enumerate}
\end{theorem}

\section{Pretty good state transfer}
%We are interested in deciding if a quantum walk can get arbitrarily close to a target state. 
A \textsl{quantum state} is a complex-valued function on the arcs of the graph represented by a unit vector. Given the transition matrix $U$ of a quantum walk and two states $x$ and $y$, we say \textsl{pretty good state transfer} occurs from $x$ to $y$ if there is a unimodular $\gamma\in\cx$ such that for any $\epsilon>0$, there is a time $t\in\ints$ with
\[\abs{U^tx -\gamma y}<\epsilon.\]
Below is an alternative definition using the inner product $\ip{U^t x}{ y}$.
\begin{lemma}\cite[Lemma 4.1]{Chan2023}\label{lem:pgstdef2}
	Given a transition matrix $U$ and two states $x$ and $y$, pretty good state transfer occurs from $x$ to $y$ if and only if for any $\epsilon>0$, there is $t\in\ints$ such that
	\[\abs{\ip{U^t x}{ y}}>1-\epsilon.\]
\end{lemma}

Since the unit circle is compact, and $\abs{\ip{U^t x}{ y}}\le 1$ for any states $x$ and $y$, the above condition is equivalent to the sequence 
\[\{\ip{U^tx}{y}: t\in \ints\}\]
having a limit point on the unit circle. We combine this observation and Theorem \ref{thm:eprojs} to simplify the analysis of pretty good state transfer from $N_t^*e_a$ to $N_t^*e_b$ in the quantum walk with respect to $(X,W)$. For ease of notation, we will abbreviate this type of state transfer as \textsl{$ab$-PGST}.

\begin{corollary}\label{cor:limpt}
	Let $W$ be a complex-weighted adjacency matrix of a connected graph $X$. Let
	\[H=\sum_{\lambda} \lambda E_{\lambda}\]
	be the spectral decomposition of the Hermitian adjacency matrix of $X$ associated with $W$. The quantum walk with respect to $(X,W)$ admits $ab$-PGST if and only if the sequence
	\[\left\{\sum_{\lambda} \cos(t\arccos \lambda) (E_{\lambda})_{ba}: t\in \ints\right \}\]
	has a limit point on the unit circle.
\end{corollary}
\proof
Let $x=N_t^*e_a$ and $y=N_t^*e_b$. By Theorem \ref{thm:eprojs},

\begin{align*}
	\ip{U^t x}{y}
	&=\sum_{\theta}e^{it\theta} \ip{F_{\theta}x}{y}\\
	&=e_b^T\left(\sum_{\theta} e^{it\theta} N_t F_{\theta}N_t^*\right) e_a\\
	&=e_b^T\left(E_1 + (-1)^t E_{-1}+(e^{i\theta}+e^{-i\theta}) \frac{1}{2} \sum_{-1<\lambda<1} E_{\lambda}\right)e_a\\
	&=\sum_{\lambda} \cos(t\arccos \lambda) (E_{\lambda})_{ba}.
\end{align*}
The result now follows from Lemma \ref{lem:pgstdef2}.
\qed

Two vertices $a$ and $b$ are \textsl{strongly cospectral relative to $H$} if for each eigenprojection $E_{\lambda}$ of $H$, there is a unimodular $\mu_{\lambda}$ such that
\begin{equation}\label{eqn:str_cosp}
	E_{\lambda} e_a = \mu_{\lambda} E_{\lambda} e_b.
\end{equation}
It was proven in \cite{Chan2023} that strong cospectrality is a neccessary condition for $ab$-PGST.

\begin{lemma}\cite[Lemma 4.3 and Lemma 6.1]{Chan2023}
	If $ab$-PGST occurs in the quantum walk with respect to $(X, W)$, then $a$ and $b$ are strongly cospectral relative to the Hermitian adjacency matrix of $X$ associated with $W$.
\end{lemma}

Following the notion in \cite{Chan2023}, we define the \textsl{eigenvalue support of vertex $a$ relative to $H$} to be the set
\[\Lambda_a=\{\lambda: E_{\lambda} e_a \ne 0\}.\]
If $a$ and $b$ are strongly cospectral, then $\Lambda_a=\Lambda_b$. Our next result shows a stronger implication of $ab$-PGST: not only does Equation \eqref{eqn:str_cosp} hold for each $\lambda$, but there is a unimodular $\gamma$ such that $\mu_{\lambda}=\pm \gamma$ whenever $\lambda\in \Lambda_a$. In this case, $\Lambda_a$ can be partitioned as
\[\Lambda_a = \Lambda_{ab}^+ \cup \Lambda_{ab}^-,\]
where
\[\Lambda_{ab}^+=\{\lambda: E_{\lambda}e_a = \gamma E_{\lambda}e_b\ne 0\}\]
and
\[\Lambda_{ab}^-=\{\lambda: E_{\lambda}e_a = -\gamma E_{\lambda}e_b\ne 0\}.\]
We will refer to $\Lambda_{ab}^+$ as the \textsl{plus set relative to $\gamma$}, and likewise $\Lambda_{ab}^-$ the \textsl{minus set relative to $\gamma$}. Note that if $H$ is real symmetric, then its eigenprojections are real, and strong cospectrality between $a$ and $b$ implies a partition of $\Lambda_a$ into plus and minus sets relative to $1$.

\begin{corollary}\label{cor:partn}
	If $ab$-PGST occurs in the quantum walk with respect to $(X,W)$, then $\Lambda_a$ can be partitioned into the plus and minus sets relative to some unimodular complex number.
%	such that
%		\[E_{\lambda} e_a =\pm \gamma E_{\lambda} e_b,\quad \forall \lambda \in \Lambda_a.\]
\end{corollary}
\proof
Since $a$ and $b$ are strongly cospectral, for each $\lambda$, there is a unimodular $\mu_r$ such that $E_{\lambda}e_a = \mu_{\lambda} E_{\lambda} e_b$. Hence
\[(E_{\lambda})_{ba} = \ip{E_{\lambda} e_a}{E_{\lambda} e_b} = \ip{\mu_{\lambda}E_{\lambda}e_b}{E_{\lambda}e_b}=\mu_{\lambda}(E_{\lambda})_{bb}.\]
Therefore
\[\sum_{\lambda} \cos(t\arccos \lambda) (E_{\lambda})_{ba}=\sum_{\lambda} \cos(t\arccos \lambda) \mu_{\lambda}(E_{\lambda})_{bb}.\]
By Corollary \ref{cor:limpt}, the sequence 
\[\left\{\sum_{\lambda} \cos(t\arccos \lambda) \mu_{\lambda}(E_{\lambda})_{bb}: t\in \ints\right\}\]
has a limit point $\gamma$ on the unit circle. It follows from
\[\sum_{\lambda} (E_{\lambda})_{bb}=1\]
that for any $\epsilon>0$, there is a $t\in\ints$ such that 
\[\abs{\cos(t\arccos\lambda)-\gamma/\mu_{\lambda}}<\epsilon,\quad \forall \lambda\in \Lambda_a.\]
As $\cos(t\arccos\lambda)$ is real, we must have $\mu_{\lambda}=\pm \gamma$ for each $\lambda\in\Lambda_a$.
\qed

The following version of Kronecker's approximation theorem is helpful in determining pretty good state transfer in coined quantum walks.

\begin{theorem}\cite{Gonek2016}\label{thm:Kronecker}
	Given $\seq{\alpha}{1}{2}{n}\in \re$ and $\seq{\beta}{1}{2}{n}\in \re$, the following are equivalent.
	\begin{enumerate}[(i)]
		\item For any $\epsilon>0$, the system
		\[\abs{q\alpha_r - \beta_r-p_r}<\epsilon, \quad r=1,2,\cdots, n\]
		has a solution $\{q, \seq{p}{1}{2}{n}\} \in \ints^{n+1}$.
		\item For any set $\{\seq{\ell}{1}{2}{n}\}$ of integers such that 
		\[\ell_1 \alpha_1 + \cdots + \ell_n \alpha_n \in \ints,\]
		we have
		\[\ell_1 \beta_1 + \cdots + \ell_n \beta_n \in \ints.\]
	\end{enumerate}
\end{theorem}

We are now ready to characterize $ab$-PGST in the quantum walk with respect to $(X,W)$ using the spectra of the associated Hermitian adjacency matrix.
\begin{theorem}\label{thm:pgst}
	Let $W$ be a complex-weighted adjacency matrix of a connected graph $X$. Let 
	\[H=\sum_{\lambda} \lambda E_{\lambda}\]
	be the spectral decomposition of the associated Hermitian adjacency matrix. The quantum walk with respect to $(X,W)$ admits $ab$-PGST if and only if the following hold.
	\begin{enumerate}[(i)]
		\item $a$ and $b$ are strongly cospectral, and their eigenvalue support can be partitioned into the plus set $\Lambda^+_{ab}$ and the minus set $\Lambda^-_{ab}$ relative to some unimodular complex number.
		\item For any set $\{\ell_{\lambda}:\lambda\in\Lambda_a\}$ of integers such that
		\[\sum_{\lambda \in \Lambda_a} \ell_{\lambda} \arccos \lambda \equiv 0 \pmod{2\pi},\]
		we have
		\[\sum_{\lambda\in\Lambda^-_{ab}}\ell_{\lambda} \equiv 0\pmod2.\]
	\end{enumerate}
\end{theorem}
\proof
Suppose first that $ab$-PGST occurs. Then by Corollary \ref{cor:partn}, there is a unimodular $\gamma$ relative to which 
\[\Lambda_a=\Lambda_{ab}^+\cup \Lambda_{ab}^-;\]
that is, for each eigenvalue $\lambda\in \Lambda_a$,  there is $\sigma_{\lambda}\in\{0,1\}$ such that
\[E_{\lambda}e_a = (-1)^{\sigma_{\lambda}} \gamma E_{\lambda}e_b.\]
Moreover, by Corollary \ref{cor:limpt}, the sequence 
\[ \left\{\gamma\sum_{\lambda} \cos(t\arccos \lambda) (-1)^{\sigma_{\lambda}}(E_{\lambda})_{bb}: t\in \ints\right\}\]
has $\gamma$ as a limit point. Thus for any $\epsilon>0$, there is a $t\in\ints$ such that 
\[\abs{\cos(t\arccos\lambda)-\cos(\sigma_{\lambda}\pi)}<\epsilon,\quad \forall \lambda\in \Lambda_a.\]
Equivalently, for any $\epsilon>0$, the system 
\[\abs{t\frac{\arccos\lambda}{2\pi}-\frac{\sigma_{\lambda}}{2}-k_{\lambda} }<\epsilon,\quad \forall \lambda\in \Lambda_a\]
has a solution $\{t\}\cup \{ k_{\lambda}:\lambda\in \Lambda_a\}\sbs \ints$. By Theorem \ref{thm:Kronecker}, this happens if and only if for any for any set set $\{\ell_{\lambda}:\lambda\in\Lambda_a\}$ of integers such that
\[\sum_{\lambda \in \Lambda_a} \ell_{\lambda} \arccos \lambda \equiv 0 \pmod{2\pi},\]
we have
\[\sum_{\lambda\in\Lambda_{ab}}\ell_{\lambda} \sigma_{\lambda}=\sum_{\lambda\in\Lambda_{ab}^-}\ell_{\lambda} \equiv 0\pmod2.\]
Thus (i) and (ii) hold. Reversing the argument shows the converse.
\qed

\section{Hypercubes}
In this section, we show every hypercube $Q_d$ admits pretty good state transfer between antipodal vertices relative to some coins constructed from a real-weighted adjacency matrix of $Q_d$. For $d$ prime, this was proven to occur with Grover coins \cite{Chan2023}.

\begin{theorem}\cite[Theorem 7.7]{Chan2023}\label{thm:prime}
	Let $p$ be a prime number. Let $A$ be the $01$-adjacency matrix of $Q_p$. The quantum walk with respect to $(Q_p, A)$ admits $ab$-PGST for any pair of antipodal vertices $a$ and $b$.
\end{theorem}

Unfortunately, this result does not generalize to all values of $d$. For example, the arc-reversal walk on $Q_4$ with Grover coins is periodic, and 
\[\max\left\{\abs{\ip{U^tN_t^*e_a}{N_t^*e_b}}: t\in \re\right\}<1.\]
Hence, a different set of coins is needed for $ab$-PGST. In the rest of this section, we construct such coins by choosing appropriate weights on the arcs of $Q_d$. Our construction exploits the fact that $Q_d$ is a Cayley graph.

Given a subset $C$ of the abelian group $\ints_2^d$, the \textsl{Cayley graph over $\ints_2^d$ with connection set $C$}, denoted $X(\ints_2^d, C)$, is the graph with vertex set $\ints_2^d$ and edge set
\[\{\{u,v\}: u, v \in \ints_2^d, u-v\in C\}.\]
For example, $X(\ints_2^d, \{e_0\})$ is isomorphic to the disjoint union of $2^{d-1}$ edges, and $X(\ints_2^d, \{\seq{e}{0}{1}{d-1}\})$ is isomorphic to the hypercube $Q_d$. It is well-known that the eigenvalues and eigenvectors of $X(\ints_q^2, C)$ are determined by the characters of $\ints_2^d$. 

\begin{lemma}\cite[Section 12.9]{Godsil1993}\label{lem:char}
	Let $X=X(\ints_2^d, C)$. For any element $g\in \ints_2$, let $\psi_g$ be the character of $\ints_2^d$ given by 
	\[\psi_g(x) = (-1)^{\ip{g}{x}}.\]
	Then $\psi_g$ is an eigenvector for $A(X)$ with eigenvalue
	\[\psi_g(C) = \sum_{c\in C} \psi_g(c).\]
	Moreover, the eigenvectors defined above are pairwise orthogonal, and they form a group isomorphic to $\ints_2^d$.
\end{lemma}

Let $m$ be a positive integer. Define two weighted adjacency matrices of $Q_d$ by
\[W_m=\sqrt{m}A\left(X(\ints_q^d, \{e_0\})\right)+\sqrt{2}\sum_{j=1}^{d-1} A(X(\ints_q^d, \{e_j\})\]
and
\[H_m = m A\left(X(\ints_q^d, \{e_0\})\right)+2\sum_{j=1}^{d-1} A\left(X(\ints_q^d, \{e_j\})\right).\]
Note that $H_m$ has constant row sum $2d-2+m$, and 
\[\frac{1}{2d-2+m} H_m\]
is the Hermitian adjacency matrix of $Q_d$ associated with $W_m$. We will find $m$ for which the quantum walk relative to $(Q_d, W_m)$ admits pretty good state transfer between the antipodal vertices. To start, we derive an explicit formula for the eigenvalues of $H_m$. Let $\wt(g)$ denote the Hamming weight of $g$.

\begin{lemma}\label{lem:Hmspectrum}
	For each element $g\in \ints_2^d$, the character $\psi_g$ given by 
	\[\psi_g(x) = (-1)^{\ip{g}{x}}\]
	is an eigenvector for $H_m$ with eigenvalue 
	\[\lambda_g = 2 d-4\wt(g) + (-1)^{g_0}(m-2).\]

\end{lemma}
\proof 
From Lemma \ref{lem:char} we see that all summands of $H_m$ are simultaneously diagonalizable by the characters of $\ints_2^d$. Hence each character $\psi_g$ is an eigenvector for $H_m$ with eigenvalue
\begin{align*}
	\lambda_g &= 2\sum_{j=0}^{d-1}\psi_g(e_j) +(m-2) \sum_{j=0}^{d-1}\psi_g(e_0)\\
	&=2\sum_{j=0}^{d-1}(-1)^{g_j} + (m-2)\sum_{j=0}^{d-1} (-1)^{g_0}\\
	&=2 d-4\wt(g) + (-1)^{g_0}(m-2). \tag*{\sqr53}
\end{align*}

Next, we establish strong cospectrality between vertices $\zero$ and $\one$.

\begin{lemma}\label{lem:parity}
Let $m$ be a positive odd integer. Let
\[H_m=\sum_{\lambda} \lambda E_{\lambda}\]
be the spectral decomposition of $H_m$. The following hold for the antipodal vertices $a=\zero$ and $b=\one$ in $Q_d$.
\begin{enumerate}[(i)]
	\item $a$ and $b$ are strongly cospectral relative to $H_m$.
%	, and their eigenvalue support can be partitioned as
%	\[\Lambda_a=\Lambda_{a b }^+ \cup \Lambda_{a b }^-,\]
%	where 
%	\[\Lambda_{a b }^+ =\{\lambda: E_{\lambda} e_{a} = E_{\lambda}e_{b}\ne 0\}\]
%	and
%	\[\Lambda_{a b }^- =\{\lambda: E_{\lambda} e_{a} =- E_{\lambda}e_{b}\ne 0\}.\]
	\item The eigenvalue support $\Lambda_a$ is symmetric about zero. Moreover, if $d$ is even, then 
	\[\lambda\in \Lambda_{ab}^{\pm}\iff -\lambda \in \Lambda_{ab}^{\pm},\]
	and if $d$ is odd, then 
	\[\lambda\in \Lambda_{ab}^{\pm} \iff -\lambda\in \Lambda_{ab}^{\mp}.\]
\end{enumerate}
\end{lemma}
\proof 
By Lemma \ref{lem:Hmspectrum}, $\lambda_g= \lambda_h$ if and only if $g_0=h_0$ and $\wt(g)=\wt(h)$. Since $\psi_g(a)=1$ and $\psi_g(b) = (-1)^{\wt(g)}$, if $\lambda_g=\lambda_h$, then we must have 
\[\frac{\psi_g(a)}{\psi_g(b)} = \frac{\psi_h(a)}{\psi_h(b)}=\pm 1.\]
Hence the projection of $e_{a}$ and $e_{b}$ onto each eigenspace of $H_m$ is either equal or opposite. This proves (i). To see (ii), let $h=\one - g$. We have 
\begin{align*}
	\lambda_h &= 2d-4\wt(h)+(m-2)(-1)^{h_0}\\
	&=2d-4(d-\wt(g))+(m-2)(-1)^{1-g_0}\\
	&=-2d+4\wt(g)-(m-2)(-1)^{g_0}\\
	&=-\lambda_g,
\end{align*}
which shows $\Lambda_a$ is symmetric about $0$. Moreover, if $d$ is even, then $\wt(g)$ and $\wt(h)$ have the same parity, from which it follows that $\psi_g(b) = \psi_h(b)$, and so $\lambda_g\in \Lambda_{ab}^{\pm}$ implies $\lambda_h\in \Lambda_{ab}^{\pm}$. Similarly, if $d$ is odd, then  $\wt(g)$ and $\wt(h)$ have different parity, and $\lambda_g\in \Lambda_{ab}^{\pm}$ implies $\lambda_h\in \Lambda_{ab}^{\mp}$.
\qed

To finish our construction, we cite two number theoretic results from \cite{Chan2023}.

\begin{lemma}\cite[Lemma 7.5]{Chan2023}\label{lem:sqff}
	Let $q,\seq{\lambda}{1}{2}{k}$ be positive integers. Suppose 
	\begin{enumerate}[(i)]
		\item $\lambda_r\in(0, q/2) \cup (d/2, q)$ for each $r$, and
		\item the square free parts of
		\[q^2-\lambda_1^2,\quad q^2-\lambda_2^2, \quad \cdots,\quad q^2-\lambda_k^2\]
		are pairwise distinct.
	\end{enumerate}
	Then the angles 
	\[\pi,\quad \arccos(\lambda_1/q),\quad \arccos(\lambda_2/q),\quad \cdots, \quad \arccos(\lambda_k/q)\]
	are linearly independent over $\rats$.
\end{lemma}

\begin{lemma}\cite[Lemma 7.6]{Chan2023}\label{lem:prime}
	Let $p$ be an odd prime. For any two distinct integers $j,k\in\{1,2,\cdots, (p-1)/2\}$, the square free parts of
	\[j(p-j),\quad k(p-k)\]
	are distinct.
\end{lemma}

We now prove our main result: every hypercube admits pretty good state transfer between antipodal vertices relative to some $W_m$.

\begin{theorem}\label{thm:cb_pgst}
	Let $d\ge 2$ be any positive integer. If $d$ is a prime, let $m=2$. If $d$ is a composite number, let $m$ be the smallest positive integer such that $2d-2+m$ is a prime. Let $a=\zero$ and $b=\one$. The quantum walk with respect to $(Q_d, W_m)$ admits $ab$-PGST.
\end{theorem}
\proof
The case for $d$ prime is settled by Theorem \ref{thm:prime}: if $m=2$, then $W_m$ is a scalar multiple of the $01$-adjacency matrix of $Q_d$. Hence we only need to prove the case where $d$ is a composite number. Let $m$ be the smallest positive integer such that $p=2d-2+m$ is a prime. Since $d\ge 2$, $p$ and $m$ must be odd, and so $a$ and $b$ are strongly cospectral vertices satisfying Theorem \ref{lem:parity} (i) and (ii). Observe from Lemma \ref{lem:Hmspectrum} that the eigenvalue support $\Lambda_a$ relative to $H_m$ is a subset of 
\[\{p-2r: r=0,1,\cdots, p\}.\]
Moreover, for each $r$, 
\[\tan\left(\arccos\left(\frac{p-2r}{p}\right)\right)\]
is a rational multiple of $\sqrt{(p-r)r}$. Hence by Lemma \ref{lem:sqff} and Lemma   \ref{lem:prime}, the angles
\begin{equation}\label{eqn:linind}
	\{\pi\}\cup \{ \arccos(\lambda/p):\lambda\in \Lambda_a, 0<\lambda<p\}
\end{equation}
are linearly independent over $\rats$. Let $\{\ell_{\lambda}:\lambda\in\Lambda_a\}$ be any set of integers such that
\begin{equation}\label{eqn:Kcond1}
	\sum_{\lambda\in \Lambda_a}\ell_{\lambda} \arccos(\lambda/p)\equiv 0\pmod{2\pi}.
\end{equation}
We aim to show 
\begin{equation}\label{eqn:Kcond2}
	\sum_{\lambda\in\Lambda_{ab}^-} \ell_{\lambda}\equiv 0 \pmod{2}.
\end{equation}

Suppose first that $d$ is even. By Lemma \ref{lem:parity}, $a$ and $b$ are strongly cospectral relative to $H_m$, and 
\[\lambda\in \Lambda_{ab}^{\pm}\iff -\lambda \in \Lambda_{ab}^{\pm}.\]
In particular, $\pm p \in \Lambda_{ab}^+$. 
Thus Equation \eqref{eqn:Kcond1} is equivalent to
\begin{align*}
	0 \equiv& \sum_{\lambda\in \Lambda_{ab}^-,\; \lambda>0} (\ell_{\lambda} \arccos(\lambda/p) + \ell_{-\lambda} \arccos(-\lambda/p))\\
	&+\sum_{\lambda\in \Lambda_{ab}^+,\; 0<\lambda<p} (\ell_{\lambda} \arccos(\lambda/p) + \ell_{-\lambda} \arccos(-\lambda/p))+\ell_{-p} \pi\\
	\equiv & \sum_{\lambda\in \Lambda_{ab}^-,\; \lambda>0} (\ell_{\lambda}- \ell_{-\lambda}) \arccos(\lambda/p) + \sum_{\lambda\in \Lambda_{ab}^+,\; 0<\lambda<p} (\ell_{\lambda}-\ell_{-\lambda}) \arccos(\lambda/p)  \\
	&+\pi \sum_{\lambda\in \Lambda_{ab},\;\lambda>0} \ell_{-\lambda} \pmod{2\pi}.
\end{align*}
Since the angles in Equation \eqref{eqn:linind} are linearly independent over $\rats$, for each  $\lambda \in \Lambda_{ab}^-$ we have $\ell_{\lambda}=\ell_{-\lambda}$. It follows that 
\[\sum_{\lambda\in \Lambda_{ab}^-} \ell_{\lambda} = \sum_{\lambda\in\Lambda_{ab}^-,\; \lambda>0}(\ell_{\lambda}+\ell_{-\lambda})\equiv 0 \pmod{2}.\]

Now suppose $d$ is odd. By Lemma \ref{lem:parity}, $a$ and $b$ are strongly cospectral relative to $H_m$, and 
\[\lambda\in \Lambda_{ab}^{\pm}\iff -\lambda \in \Lambda_{ab}^{\mp}.\]
In particular, $p\in \Lambda_{ab}^+$ and $-p\in \Lambda_{ab}^-$. Thus Equation \eqref{eqn:Kcond1} is equivalent to
\begin{align*}
	0 \equiv& \sum_{\lambda\in \Lambda_{ab}^+\backslash\{p\}} \ell_{\lambda} \arccos(\lambda/p) + \sum_{\lambda\in \Lambda_{ab}^+\backslash\{p\}} \ell_{-\lambda} \arccos(-\lambda/p)+\ell_{-p}\pi \\
	\equiv & \sum_{\lambda\in \Lambda_{ab}^+\backslash\{p\}} (\ell_{\lambda}-\ell_{-\lambda}) \arccos(\lambda/p) + \pi \sum_{\lambda\in \Lambda_{ab}^+} \ell_{-\lambda}  \pmod{2\pi}.
\end{align*}
Since the angles in Equation \eqref{eqn:linind} are linearly independent over $\rats$, for each $\lambda\in\Lambda_{ab}^+\backslash\{p\}$ we have $\ell_{\lambda}=\ell_{-\lambda}$, and 
\[\sum_{\lambda\in \Lambda_{ab}^+} \ell_{-\lambda}  \equiv 0 \pmod{2}.\]
It follows that
\[\sum_{\lambda\in \Lambda_{ab}^-} \ell_{\lambda} = \sum_{\lambda\in \Lambda_{ab}^+} \ell_{-\lambda}  \equiv 0 \pmod{2}.\]

In both cases, Equation \eqref{eqn:Kcond2} holds, and so by Theorem \ref{thm:pgst}, the quantum walk with respect to $(Q_d, W_m)$ admits $ab$-PGST.
\qed

\section{A sufficient condition for $ab$-PGST}
Our approach suggests a sufficient condition for other graphs to admit pretty good state transfer with real coins.

\begin{theorem}\label{thm:general}
	Let $X$ be a connected graph. Let $H$ be a symmetric non-negative adjacency matrix of $X$, with spectral radius $p$ for some prime $p$. Let $W=H^{\circ 1/2}$. Suppose $a$ and $b$ are strongly cospectral relative to $H$, and $\Lambda_a$ is a subset of 
	\[\{p-2r: r=0,1,\cdots,p\}.\]
	Then $ab$-PGST occurs in the quantum walk with respect to $(X,W)$ if one of the following holds.
	\begin{enumerate}[(i)]
		\item For any pair $\lambda, -\lambda$ in $\Lambda_a$, 
		\[\lambda\in \Lambda_{ab}^{\pm} \iff -\lambda \in \Lambda_{ab}^{\pm}.\]
		\item For any pair $\lambda, -\lambda$ in $\Lambda_a$, 
		\[\lambda\in \Lambda_{ab}^{\pm} \iff -\lambda \in \Lambda_{ab}^{\mp}.\]
	\end{enumerate}
\end{theorem}
\proof
Since $H$ is non-negative, by the Perron–Frobenius theorem, $p\in \Lambda_{ab}^+$. If $-p\in \Lambda_{ab}^-$, then $X$ is bipartite, and an similar argument to Theorem \ref{thm:cb_pgst} shows $ab$-PGST. Otherwise, the condition
\[\sum_{\lambda\in\Lambda_a} \lambda \arccos(\lambda/p)\equiv 0\pmod{2\pi}\]
together with linear independence of the angles 
\[	\{\pi\}\cup \{ \arccos(\lambda/p):\lambda\in \Lambda_a, 0<\lambda<p\}\]
over $\rats$ imply that
\begin{enumerate}[(a)]
	\item $\ell_{\lambda}=\ell_{-\lambda}$ for any pair $\lambda, -\lambda\in \Lambda_a$.
	\item $\ell_\lambda=0$ for any $\lambda\in\Lambda_a$ such that $-\lambda\notin \Lambda_a$.
	\item $\displaystyle \sum_{\lambda: \lambda\in\Lambda_{ab} \text{ and } -\lambda\in \Lambda_{ab}} \ell_{\lambda}$ is even.
\end{enumerate}
Thus, if (i) or (ii) holds, then 
\[\sum_{\lambda\in\Lambda_{ab}^-} \ell_{\lambda}\equiv 0\pmod{2}.\]
By Theorem \ref{thm:pgst}, the quantum walk with respect to $(X,W)$ admits $ab$-PGST.
\qed

\section{Future work}
We have found real quantum coins that enable $ab$-PGST on every hypercube. Our coins are slight modifications of Grover coins that require change of the weight on one arc per vertex. As indicated by Theorem \ref{thm:general}, it is likely that a similar construction works for a much broader class of graphs. An interesting direction is to determine, for each family of graphs that do not admit pretty good state transfer with Grover coins, the minimum change of weights needed for this phenomenon to occur.

\section*{Acknowledgement}
This material is based upon work supported by the National Science Foundation under Grant No. 2348399.

\bibliographystyle{amsplain}
\bibliography{qw.bib}

\end{document}